\documentclass[12pt,a4paper]{article}

\usepackage{amsmath,amsthm,mathrsfs,amsfonts}
\usepackage{ulem}	
\usepackage{hyperref}
\usepackage{amsfonts,todonotes}
\usepackage{amssymb}
\usepackage{bbm}
\usepackage{xcolor}

\usepackage{soul}
\usepackage{graphicx}
\usepackage{subfig}
\usepackage{float}

\def\defeq{\mathrel{\mathop:}=}

\newcommand{\I}{\mathbbm{1}}

\usepackage{multirow}
\newcommand*{\QEDB}{\hfill\ensuremath{\square}}%

\usepackage{parskip}
\newtheorem{theorem}{Theorem}
\newtheorem{lemma}{Lemma}
\newtheorem{definition}{Definition}

\newtheorem{corollary}{Corollary}

\usepackage{booktabs}
\usepackage[shortlabels]{enumitem}

\begin{document}

	\title{\bf Level sets of  depth measures and central dispersion in abstract spaces }
	
	\author{ Alejandro Cholaquidis \thanks{acholaquidis@cmat.edu.uy} \hspace{0.1cm} \, \hspace{.2cm}   Ricardo Fraiman \thanks{rfraiman@cmat.edu.uy}  \hspace{.2cm}    Leonardo Moreno \thanks{mrleo@iesta.edu.uy}
	}
	\maketitle

\begin{abstract}
The lens depth of a point has been recently extended to general metric spaces, which is not the case for most depths.
It is defined as the probability of being included in the intersection of two random balls centred at two random points $X$  and $Y$, with the same radius $d(X,Y)$. We study the consistency in Hausdorff and measure distance, of the level sets of the empirical lens depth, based on an iid sample  on a general metric space.
We also prove that the boundary of the empirical level sets are consistent estimators of their population counterparts, and analyze two real-life examples.
\end{abstract}

\section{Introduction}
 Statistical depth functions have gained  importance in the last three decades. The best well-known and most studied among them  include the half space depth \cite{tukey1975}, the simplicial depth \cite{liu1990,liu1992}, and the multivariate $L^1$ depth. 
 
Other well-known depths are the convex hull peeling depth \cite{barnett1976}, the  Oja depth \cite{oja1983}, and the spherical depth \cite{elmore2006}, among others.
They have been extended to functional spaces, see \cite{fraiman2001,lopez2009,claeskens2014,cuevas2009}, to  Riemannian manifolds, see \cite{fraiman2019}, and also to general metric spaces, see \cite{fraiman2020}.
 Several different applications of depths notions have been proposed, in particular for classification problems, by means of the depth-depth method \cite{vencalek2017}, or to functional data, see \cite{mosler2017}.
   Most classical notions of depth, introduced initially on $\mathbb{R}^d$, can not be directly extended to general metric spaces or even to functional spaces or manifolds. \color{black}
Some of them are computationally infeasible on high dimensional spaces, such as Liu's or Tukey's depth, because the computational complexity is exponential in the number of dimensions.
 This is not the case of the depth introduced in \cite{liu2011}, i.e. the lens depth, whose computational complexity is of order $n^2$, and,   as we will see, \color{black} it can be easily extended to general metric spaces.
This makes the lens depth particularly suitable for estimating it's level sets, by means of the level sets of it's empirical version, based on an iid sample, which is one of the main goals of this manuscript.
   An extension to Riemannian manifolds of the lens depth (called weighed lens depth)  was recently introduced in \cite{cholaquidis2020} to tackle the supervised classification problem, and a point-wise a.s consistency result is obtained (see Theorem 2). Here we focus on a different problem, i.e: level set estimation of the lens depth on general metric spaces. This require uniform a.s. consistency of the empirical lens depth to the population one, which is developed in Section \ref{cons}.\color{black}
   
Level set estimation of depths was initially studied in \cite{tukey1975}, as a key tool for the visualization and exploration of data.  Other significant contributions can be found in  \cite{koshevoy1997, serfling2000, serfling2002b, dyckerhoff2016}).  As it was pointed out in \cite{liu1999}, ``the shape and size of these  levels, as well as  the direction and speed at which they expand, provides insight into the dispersion, kurtosis, and asymmetry of the underlying distribution".
It also allows to extend the notion of quantiles to multivariate or functional data, and can be used for outlier detection, see for instance  \cite{febrero2008, dai2019}, as well as for supervised classification, see \cite{ruts1996,hubert2017}. 


More formally, let $X_1,X_2$ be two random variables defined on a (rich enough) probability space $(\Omega, \mathcal{A},\mathbb{P})$, taking  values in a complete separable metric space $(M,d)$ endowed with the Borel $\sigma$-algebra.
Assume that they are independent and identically distributed.   In what follows,  the distribution of a random element $X$ of $M$ will be denoted by $P_X$. 

\begin{figure}[!ht] 
	\centering
	\subfloat{\includegraphics[width=70mm]{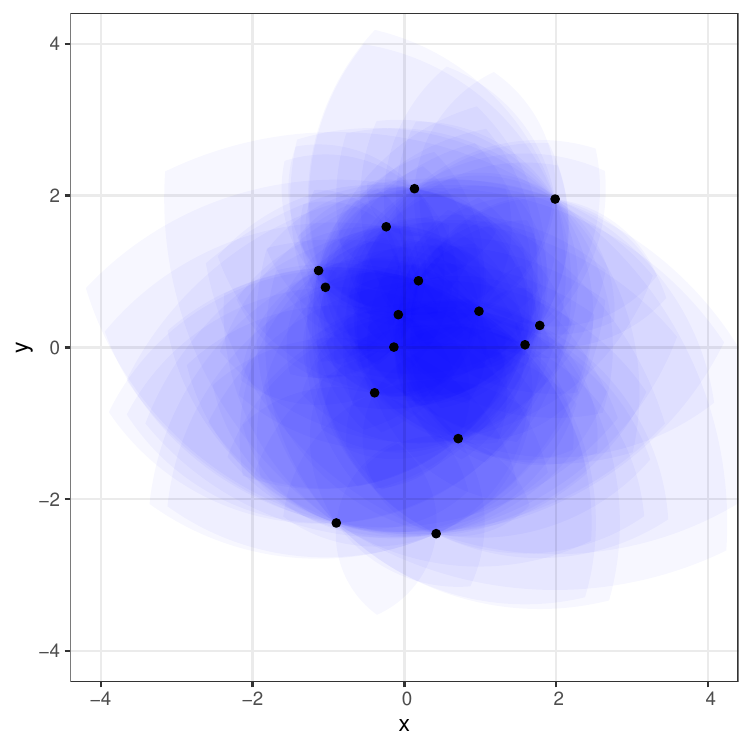}}
	\caption{Level sets of the lens depth based on a sample of 15 iid random vectors, distributed as a bivariate standard normal distribution.
	     The intensity of blue represents the depth.} \label{ldc}  
\end{figure}

Given $x_1, x_2 \in M$ define their ``associate lens" by
$$A(x_1, x_2):= B(x_1, d(x_1,x_2)) \cap B(x_2, d(x_1,x_2)),$$ where $B(p,r)$ is the closed ball centred at $p$ with radius $r>0$. The lens depth of a point $x\in M$ is defined by  $\textrm{LD}(x)=\mathbb{P}(x\in A(X_1,X_2))$. \color{black}
Given an iid sample $\mathcal{X}_n= \{X_1, \ldots, X_n\}$  from a distribution  $P_X$, the empirical version of $\textrm{LD}$ is given by the $U$-statistics of order two,  
\begin{equation}
	\widehat{\textrm{LD}}_n(x)=\binom{n}{2}^{-1} \sum_{1 \leq i_1 < i_2 \leq n} \I_{A(X_{i_1}, X_{i_2})}(x).
	\label{estimadorprofundidad_c}
\end{equation}
  To gain some insight into the shape of $\{\widehat{\textrm{LD}}_n \geq \lambda\}$, see Figure \ref{ldc}. \color{black} 

The general approach to obtain consistency results (w.r.t.\ Hausdorff distance) for the plug-in estimator $\{\widehat{\textrm{LD}}_n \geq \lambda\}$ to its population counterpart $\{\textrm{LD}\geq \lambda\}$ is  to prove the uniform convergence of $\widehat{\textrm{LD}}_n$ to $\textrm{LD}$.
  
In general metric spaces it will be necessary to restrict the convergence to compact sets.
This is proved in Theorem \ref{convunif}, but the results also hold on $\mathbb{R}^d$ under milder conditions, see Theorem \ref{Rd}.

This paper is organized as follows.  Section \ref{notation} introduces the notation, some previous definitions,   and states two important results given in Theorems \ref{molch} and \ref{tmcu2006} that will be used to prove our main results. \color{black}
The almost surely (a.s.) uniform consistency of $ \widehat{\textrm{LD}}_n$ is stated in Section \ref{cons}, while its asymptotic distribution is formulated in Section \ref{asymptotic}.
The a.s. consistency in Hausdorff and measure distance,  for the empirical level sets $\{\widehat{\textrm{LD}}_n\geq \lambda\}$, as well as the a.s. consistency of its boundary, is stated in Section \ref{levsec}.
Lastly, in Section \ref{realdata} we  tackle the study of two interesting real datasets,  the vectocardiogram dataset (see subsection \ref{vectcard}) and the influenza dataset (see subsection \ref{influ}).
All proofs are given in the Appendix.

\section{Preliminaries}\label{notation}

In this section we will introduce the notation and necessary definitions  used throughout this paper.
Given a metric space $(M,d)$, which will be assumed to be separable, complete and locally connected, we denote by $B(x,r)$ the closed ball centred at $x$ with radius $r>0$.
The boundary of a set $A\subset M$ is denoted by $\partial M$.
Given $x\in A^c$, $d(x,A)=\inf_{a\in A} d(x,a)$ and $\textrm{diam}(A)=\sup_{x,y\in A}d(x,y)$.
Given two closed sets $A,B$, the Hausdorff distance between them is defined as   
\begin{equation}
d_H(A,B) \defeq \max\Big\{\max_{a\in A}d(a,B), \ \max_{b\in B}d(b,A)\Big\}
\end{equation}
Given a function $f: M \rightarrow \mathbb{R}$ and $\lambda \in \mathbb{R}$, we denote by $\{f\geq \lambda\}$ the $\lambda$-level set $\{x\in M: f(x)\geq \lambda\}$.

For the estimation of level sets in general metric spaces, two main results will play a key role.
The first one is Theorem 2.1 in \cite{molchanov1998}. We will make use of the following slightly restricted version.

\begin{theorem}[Molchanov, (1998)]
	\label{molch}	
Let $f_n,f:M\to \mathbb{R}$ be continuous functions. Assume that for each compact set $K_0$, $\sup_{x\in K_0}|f_n(x)-f(x)|\to 0$.
Assume that for all $\lambda\in [c_1,c_2]$, $\{f\geq\lambda\}\subset \overline{\{f>\lambda\}}$.
 Then 
$$\sup_{c_1\leq \lambda \leq c_2} d_H\big(\{f_n\geq \lambda \}\cap K_0,\{f\geq \lambda\}\cap K_0\big)\to 0.$$
\end{theorem}
The convergence of the boundaries of the level sets is in general more involved. To prove that, we will use the following result.
\begin{theorem}[Cuevas, Goz\'alez-Manteiga, Rodr\'iguez-Casal, (2006)]
\label{tmcu2006}
Given a continuous function  $f:M\to \mathbb{R}$,  let $(\Omega,\mathcal{A},P)$ be a probability space and $f_n=f_n(\omega,\cdot)$, with $\omega\in \Omega$, a sequence of random functions, $f_n:M\to \mathbb{R}$, $n=1,2,\dots,$.
Assume that for each $n$, $f_n$ is continuous with probability one.
Assume that the following assumptions are fulfilled.

\begin{itemize}
\item[$h_1)$]  $M$ is locally connected.
\item[$h_2)$] For all $x\in \partial \{f\geq \lambda\}$, there exist sequences $u_n,l_n\to x$ such that
$f(u_n)>\lambda$ and $f(l_n)< \lambda$.
\item[$f_1)$] $\partial\{ f \geq \lambda \}\neq \emptyset$.
Moreover, there exists a $\lambda^{-}< \lambda$ such that the set $\{ f \geq \lambda^{-} \}$ is compact.
\end{itemize}
If $\sup_{ x \in M} \big \vert f(x) - f_n(x) \big \vert \rightarrow 0$ a.s., then 

$$d_H \big(\partial\{f\geq \lambda\} , \partial \{f_n\geq \lambda\} \big) \to 0, \quad a.s., \textrm{ as } n \rightarrow \infty.$$

\end{theorem}

\section{Uniform consistency  of $\widehat{\textrm{LD}}_n$} \label{cons}

As mentioned above, the key points to prove the consistency in Hausdorff distance of the level sets of $\widehat{\textrm{LD}}_n$ are Theorems \ref{tmcu2006} and \ref{molch}.
Then we have to prove that $\widehat{\textrm{LD}}_n$ converges uniformly to $\textrm{LD}$ a.s., which is the main goal of this section.

To obtain the a.s. uniform convergence of $\widehat{\textrm{LD}}_n$, we will use the following version of Theorem 1 in \cite{billingsley1967}.

\begin{theorem}[Billingsley and Tops{\o}e, (1967)]
	\label{bill}	
Let $\mathcal{B}(M\times M)$ be the class of all real valued, bounded, measurable functions defined on the metric space $(M\times M,\rho)$, where $\rho(z,y)=\max\{d(z_1,y_1),d(z_2,y_2)\}$.
Suppose $\mathcal{F}\subset \mathcal{B}(M\times M)$ is a subclass of functions.
Then 
\begin{equation*}
\sup_{f\in \mathcal{F}} \Big| \int fdP_n-\int fdP\Big|\rightarrow 0,
\end{equation*}
for every sequence $P_n$ that converges weakly to $P$ if, and only if,
$$\sup\{|f(z)-f(t)|:f\in \mathcal{F},z=(z_1,z_2),t=(t_1,t_2)\in M\times M\}< \infty,$$
and for all $\epsilon>0$,
\begin{equation}\label{biltop}
\lim_{\delta\rightarrow 0} \sup_{f\in \mathcal{F}} P\Big[\{y=(y_1,y_2): \omega_f\{B_\rho(y,\delta)\}\geq \epsilon\}\Big]=0,
\end{equation}
where $\omega_f(A)=\sup\{|f(z)-f(t)|:z,t\in A\}$ and $B_\rho(y,\delta)$ is the open ball in the metric space $(M\times M,\rho)$ of radius $\delta>0$.
\end{theorem}

To prove that \eqref{biltop} is fulfilled, we will use the following lemma.

\begin{lemma}\label{lemaux}
	Let $x\in M$ and  $y=(y_1,y_2)\in M\times M$ be such that $d(x,\partial A(y_1,y_2))>3\delta$. Then $\omega_{f_x}\{B_\rho(y,\delta)\}=0$, where $f_x(t_1,t_2)=\mathbb{I}_{A(t_1,t_2)}(x)$ with $t_1,t_2\in M$.
\end{lemma}

\begin{theorem} \label{convunif} Let $(M,d)$ be a complete separable metric space and $P_X$ a Borel measure on $M$. Assume that $P_X(\partial B(x,\delta))=0$ for all $x\in M$ and $\delta>0$.
For any compact set $K\subset M$,
	$$\sup_{x\in K} |\widehat{\textrm{LD}}_n(x)-\textrm{LD}(x)|\to 0\quad a.s., \text{ as } n\to \infty.$$
\end{theorem}

As a direct consequence of the previous results, if $M$ is a compact manifold, we have that $\sup_{x\in M} |\widehat{\textrm{LD}}_n(x)-\textrm{LD}(x)|\to 0$.
For $\mathbb{R}^d$, uniform convergence can be obtained, as stated in the following theorem.

\begin{theorem}\label{Rd} Assume that $P_X(\partial B(x,\delta))=0$ for all $x\in \mathbb{R}^d$ and $\delta>0$.
 Then
	 $$\sup_{x\in \mathbb{R}^d} |\widehat{\textrm{LD}}_n(x)-\textrm{LD}(x)|\to 0\quad a.s., \text{ as } n\to \infty.$$
\end{theorem}

\section{Asymptotic distribution}\label{asymptotic}

To obtain the asymptotic law of $\widehat{\textrm{LD}}_n$ we will use Proposition 10 of \cite{gine1996}, which is a very general result for empirical processes, applied to $U$ statistics.
  This requires proving that the family of sets $\{A(x,y), x, y \in M\}$ has finite Vapnik--Chervonenkis (VC) dimension (see \cite{devroye2013}). \color{black}
To prove the result, we introduce the  \textit{Vapnik-ball} condition.

\begin{definition} A metric space $(M,d)$ fulfills  \textit{Vapnik-ball} condition if the the family $\{F_p\}_{p\in M}$ where $F_p=\{(x,y) \in M \times M : d(p,x) < d(x,y) \}$ has finite VC dimension varying $p$.
\end{definition} 

Let us introduce some notation and general ideas of empirical processes, following \cite{arcones1993}.
To each $x\in M$ we associate the function $f_x:M^2\to \mathbb{R}$ defined by $f_x(a,b)=\mathbb{I}_{A(a,b)}(x)$.
The family of all these functions is denoted by $\mathcal{F}$.
Then we can define $\widehat{\textrm{LD}}_n$, a function of the random variable $X$ and the point $x$, by 
$$\widehat{\textrm{\textrm{LD}}}_n(P_X,f_x)=\binom{n}{2}^{-1}\sum_{1\leq i<j\leq n} f_x(X_i,X_j),$$
where $X_1,\dots, X_n$ is an iid sample from $X$ with distribution $P_X$.
Denote the metric space of all real valued bounded functions defined on $\mathcal{F}$, endowed with the supremum norm, by $l^\infty(\mathcal{F})$.
We can consider $\widehat{\textrm{LD}}_n$ as a random variable with values in $l^\infty(\mathcal{F})$ given by $\omega\in \Omega\to \widehat{\textrm{LD}}_n:\mathcal{F}\to \mathbb{R}$.
We define $\textrm{LD}\in l^\infty(\mathcal{F})$ by $\textrm{LD}(f_x)=\mathbb{P}(f_x(X_1,X_2))$, with $X_1,X_2$ being independent copies of $X$.
We want to derive the limit law of the random element of $l^\infty(\mathcal{F})$ given by $G_n=\sqrt{n}(\widehat{\textrm{LD}}_n-\textrm{LD})$.
The limit law of $G_n$ is denoted by $2GP_X$, where $GP_X$ is a Brownian bridge associated to $P_X$,   which is a random element  in $l^\infty(\mathcal{F})$, and is defined through its finite dimensional distributions.
More precisely, for all $k$, $(GP_X(f_{x_1}),\dots,GP_X(f_{x_k}))$   is a zero-mean Gaussian process whose $k\times k$ covariance matrix  has $(i,j)$th element
\begin{equation}\label{cov}
E(GP_X(f_{x_i})GP_X(f_{x_j}))= P_2(f_{x_i}f_{x_j})- P_2(f_{x_i})P_2(f_{x_j}),
\end{equation}
where $$P_2(f_x) \defeq \int_{M^2} f_x(x_1,x_2) dP_X(x_1) \times dP_X(x_2).$$
The convergence of $G_n$ to $2GP_X$ is in law, which in this case is equivalent to the convergence in distribution of the random vector $(G_n(f_{x_1}),\dots,G_n(f_{x_k}))$ to $N(0,\Sigma)$, where the $k\times k$ matrix $\Sigma$ has $(i,j)$th element given by \eqref{cov}.

	\begin{theorem} {(Limit Law)}
		\label{prop3} Let $M$ be a compact manifold fulfilling the  \textit{Vapnik-ball} condition.
Following the previous notation, the stochastic processes
		 $G_n$ converge in law to $2GP_X$, as  $n\to \infty$.
	\end{theorem}  

\section{Level set estimation}\label{levsec}

Given positive numbers $c_1<c_2$,  assume that for all $\lambda\in [c_1,c_2]$, we have $\{f\geq\lambda\}\subset \overline{\{f>\lambda\}}$.
From Theorem 2.1 in \cite{molchanov1998} together with Theorem \ref{convunif}, and the fact that $\textrm{LD}$ is a continuous function, we get
$$\sup_{c_1\leq \lambda\leq c_2} d_H\big(\{\widehat{\textrm{LD}}_n\geq \lambda\}\cap K,\{\textrm{LD}\geq \lambda\}\cap K\big)\to 0\quad a.s., \text{ as } \quad n\to \infty$$  
for any compact set $K$.
If $\nu(\{\textrm{LD}=\lambda\})=0$, it follows easily from Theorem \ref{convunif} that for all compact $K$,
$$d_\nu\big(\{\widehat{\textrm{LD}}_n\geq \lambda\}\cap K,\{\textrm{LD}\geq \lambda\}\cap K\big)\to 0\quad a.s., \text{ as } \quad n\to \infty.$$
The following theorem states that $\partial\{\widehat{\textrm{LD}}_n\geq \lambda\} \cap K$ is a consistent estimator of $\partial \{\textrm{LD}\geq \lambda\}\cap K$.
The proof follows the same lines used to prove Theorem 1 in \cite{cuevas2006}.
However, we can not apply that theorem directly because $\widehat{\textrm{LD}}_n(x)$ is not a continuous function, since the range of $\widehat{\textrm{LD}}_n(x)$ is contained in the set $\{k\binom{n}{2}^{-1}:k=0,\dots,\binom{n}{2}\}$.

\begin{theorem} \label{levels} Let $\lambda>0$ be such that $\{ \textrm{LD}\geq\lambda\}\neq \emptyset$. Under the assumptions of Theorem \ref{convunif}, together with hypotheses $h1$ and $h2$ of Theorem \ref{tmcu2006}, we have that

	\begin{equation*}\label{hbdr}
		\lim_{n\to \infty} d_H\big(\partial\{\widehat{\textrm{LD}}_n\geq \lambda\} \cap K,\partial \{\textrm{LD} \geq \lambda\}\cap K\big)=0\quad a.s.,
		\end{equation*}
	for all compact sets $K$.
\end{theorem}

For the metric space $(\mathbb{R}^d,\|\cdot\|)$ endowed with the Euclidean norm we have the following trivial corollary, 

\begin{corollary} Assume that hypothesis h2 of Theorem \ref{tmcu2006} is fulfilled.
Assume also that $\{\textrm{LD}\geq \lambda^{-}\}$ is non-empty and compact, for some $\lambda^{-}<\lambda$. Then
		\begin{equation}\label{hbdr2}
		\lim_{n\to \infty} d_H(\partial\{\widehat{\textrm{LD}}_n\geq \lambda\},\partial \{\textrm{LD} \geq \lambda\})=0\quad a.s.
	\end{equation}
	
\end{corollary}

\section{An application to the study of two sets of real data}\label{realdata}
In what follows we analize two interesting real data sets: the vectocardiogram and the influenza data sets. \color{black}

\subsection{The vectocardiogram dataset}\label{vectcard}

\begin{figure}[!ht] 
	\centering
	\subfloat{\includegraphics[width=80mm]{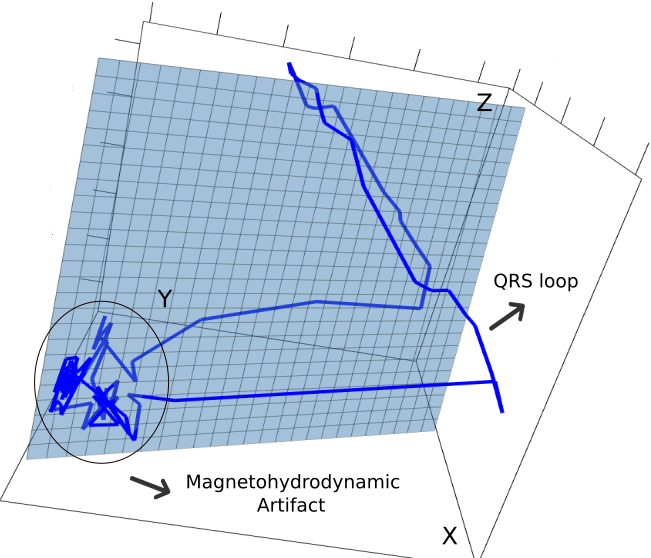}}
	\caption{QRS curve loop.
		The plane in blue is the best two dimensional approximation of the curve by least squares.}
	\label{dibujo}  
\end{figure}

We consider a real life dataset, where the data belong to the Stiefel manifold $SO(3,2)$ of all orthonormal $2$-frames in $\mathbb{R}^3$ considered as $3\times 2$ orthogonal matrices (see \cite{hatcher}).
The dataset consists of 98 vectocardiograms from children  with ages varying between 2 and 19.
Vectocardiography is a method that produces a three dimensional curve which comprises  the  records of the magnitude and direction of the electrical forces generated by the heart over time.
These curves are called  QRS loops,  see Figure \ref{dibujo}.
In \cite{downs1969} there is associated to each curve an element of $SO(3,2)$ that represents some of the information of the curve, see also \cite{downs1972}.

This sample has been previously analysed in the literature, see for instance  \cite{chikuse2012, chakraborty2019, pal2019}.
A very important problem is outlier detection, corresponding to children with (possibly) cardiological problems.

\begin{figure}[!ht] 
	\centering
	\subfloat{\includegraphics[width=110mm]{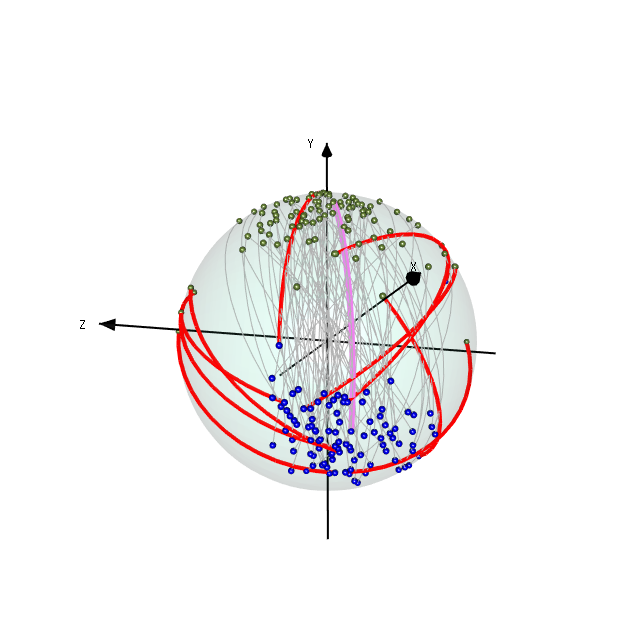}}
	\caption{Data visualization of the data from  $SO(3,2)$.
		The violet arc represents the deepest observation according to  lens depth.
		The outliers represented by red arcs correspond to the data outside the   $\lambda=0.10$ level set of the depth.
	} \label{bola}  
\end{figure}

Figure \ref{bola} represents each matrix in $SO(3,2)$ as two points in $S^2$, one for each column.
They are joined by an arc in $S^2$.
The arc joining the deepest pair of observations (w.r.t. lens depth) is represented in violet, while the outliers (for a level $\lambda=0.10$) are represented as red arcs.

\subsection{Influenza dataset}\label{influ}

This dataset consists of longitudinal data of the influenza virus, belonging to the family Orthomyxoviridae.
  An important problem analized in
 \cite{fraiman2020},  is  to model the genomic evolution of the virus, see \cite{smith2004mapping}. 
In this paper we focus on another important issue: modelling the temporal variability of the virus by means of the lens depth, and use this to predict and anticipate a possible pandemic. \color{black}
The influenza virus has an RNA genomic which is very common: it produces diseases like yellow fever and hepatitis and annually costs half a million deaths worldwide.
It is well known that these viruses change their genetic pattern over time, which is vital for developing a possible vaccine.
We will study the the H3N2 variant of the virus, in particular, the subtype hemagglutini (HA), which produced the SARS pandemic in 2002.
This variant  is known to have a variability in its genetic arrangement over  time, see \cite{altman2006, monod2018}.

The dataset can be found in GI-SAID $\textrm{EpiFlu}^{TM}$ database \footnote{\url{www.gisaid.org}}, providing 1089 genomic sequences of H3N1 from 1993 to 2017 in New York, aligned using MUSCULE, see  \cite{edgar2004}.
We used reduced trees of 5 leaves, as in \cite{monod2018}, to capture the structure of the data.
This set of trees can be endowed with a distance,
see \cite{billera2001}.
The database used was obtained from the  GitHub repository  \url{https://github.com/antheamonod/FluPCA}.

For this kind of data, constructing measures of centrality and variability, as well as confidence regions, is a problem that has been previously addressed in the literature, see for instance  \cite{barden2018, brown2020, willis2019} and references therein.

\begin{figure}[!ht] 
	\centering
	\subfloat{\includegraphics[width=85mm]{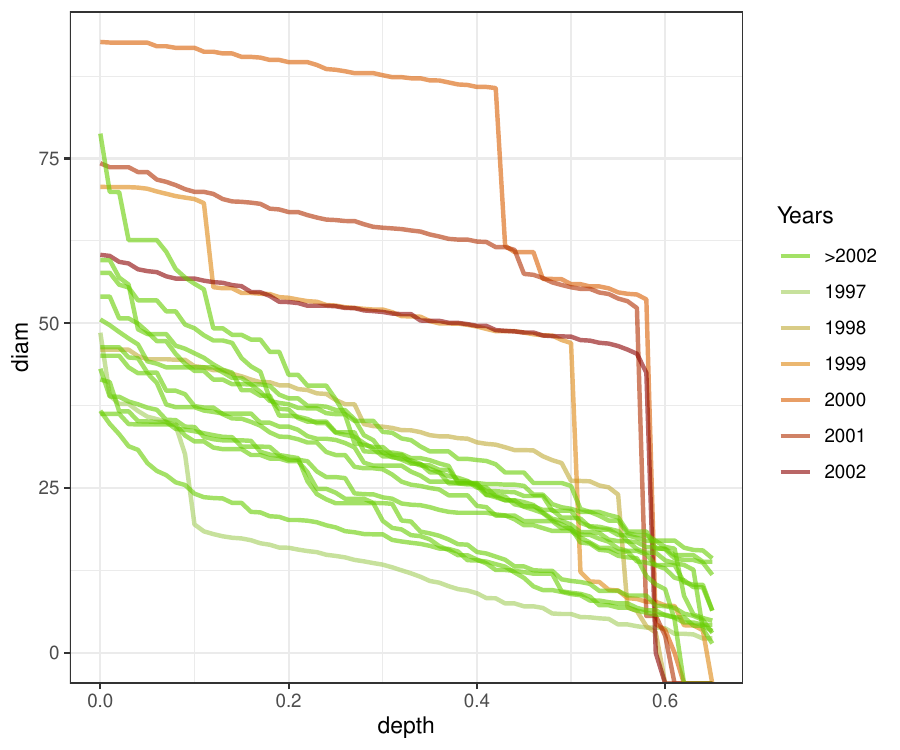}}
	\caption{Diameters of the lens depth level sets from the sample of trees for the different years, with respect to $\lambda$.} \label{disp}  
\end{figure}

For each year we computed the empirical lens depth of the trees, considered on the manifold of phylogenetic trees.
We estimated the diameter of the level sets from the sample points that belong to the level sets, see Figure \ref{disp}.
As can be seen, there is a larger dispersion  in the years prior to the pandemic of 2002 compared to later years.

\section*{Appendix}

\textit{Proof of Lemma \ref{lemaux}}

If $z=(z_1,z_2)\in  B_\rho(y,\delta)$, then $d(z_1,y_1)<\delta$ and $d(z_2,y_2)<\delta$.
Let us prove first that
	\begin{equation}\label{inc1}
		\{x\in A(y_1,y_2):d(x,\partial A(y_1,y_2))>2\delta\}\subset A(z_1,z_2)
	\end{equation}
	for all $z=(z_1,z_2)\in  B_\rho(y,\delta)$.
If $x\in A(y_1,y_2)$ and $d(x,\partial A(y_1,y_2))>2\delta$, then $d(x,y_1)\leq d(y_1,y_2)-2\delta\leq d(y_1,z_1)+d(z_1,z_2)+d(z_2,y_2)-2\delta\leq  d(z_1,z_2)$.
In the same way, if $d(x,y_2)\leq d(z_1,z_2)$, then $x\in A(z_1,z_2)$, which implies \eqref{inc1}.
Next let us prove that 
	\begin{equation}\label{inc2}
		A(z_1,z_2)\subset \{x: d(x,A(y_1,y_2))\leq 3\delta\}
	\end{equation}
	For any $x\in A(z_1,z_2)$, $d(x,y_1)\leq d(x,z_1)+d(z_1,y_1)\leq d(z_1,z_2)+\delta\leq d(y_1,y_2)+3\delta$.
In the same way,
	$d(x,y_2)\leq d(x,z_2)+\delta \leq d(z_1,z_2)+\delta\leq d(y_1,y_2)+3\delta$.
Then $x\in  \{x: d(x,A(y_1,y_2))\leq 3\delta\}$.
To prove that $\omega_{f_x}\{B_\rho(y,\delta)\}=0$, it is enough to prove that for all  $(z_1,z_2),(t_1,t_2)\in B_\rho((y_1,y_2),\delta)$
	\begin{equation}\label{eq2}
		|\mathbb{I}_{A(z_1,z_2)}(x)-\mathbb{I}_{A(t_1,t_2)}(x)|=0.
	\end{equation}
Observe that  if  $x\in A(y_1,y_2)$ and $d(x,\partial A(y_1,y_2))>3\delta$, then by \eqref{inc1}, $x\in A(z_1,z_2)$, so $\mathbb{I}_{A(z_1,z_2)}(x)=0$ and $x\in A(t_1,t_2)$ and hence $\mathbb{I}_{A(t_1,t_2)}(x)=0$ which then shows that \eqref{eq2} holds.
Proceeding in the same way if $x\notin A(y_1,y_2)$ and $d(x,\partial A)>3\delta$, by \eqref{inc2} $x\notin A(z_1,z_2)$, which implies that $\mathbb{I}_{A(z_1,z_2)}(x)=0$.
Also $x\notin A(t_1,t_2)$, which implies that $\mathbb{I}_{A(t_1,t_2)}(x)=0$, so again \eqref{eq2} holds.
\QEDB

\textit{Proof of Theorem \ref{convunif}}

	We will apply Billingsley's theorem to the set of functions $$\mathcal{F}=\{f_x(z,y)=\mathbb{I}_{A(z,y)}(x):x\in K\},$$ 
	where the sequence $P_n$ of probability measures on $M\times M$ is such that $P_n(X_i,X_j)=(1/2)\binom{n}{2}^{-1}$ if $i\neq j$ and $P_n(X_i,X_j)=0$ if $i=j$.
Let $P$ be the product measure $P_X\times P_X$.
In this case
	$$\sup_{f\in \mathcal{F}} \Big| \int fdP_n-\int fdP\Big|=\sup_{x\in K} |\widehat{\textrm{LD}}_n(x)-\textrm{LD}(x)|.$$
	
	Clearly,  $\sup\{|f(z)-f(t)|:f\in \mathcal{F},z,t\in M\times M\}=2$.
So we have to prove that \eqref{biltop} holds. By  Lemma \eqref{lemaux}, to prove  \eqref{biltop}, it is enough to prove that 
	$$\lim_{\delta\to 0} \sup_{x\in K} P_X\times P_X\{(y_1,y_2): d(x,\partial A(y_1,y_2))\leq 3\delta\}=0.$$
	By the dominated convergence theorem, since $P_X(\partial B(x,\delta))=0$ for all $\delta>0$ and all $x\in M$, it follows that for fixed $\delta$, $P_X\times P_X\{(y_1,y_2): d(x,\partial A(y_1,y_2))\leq 3\delta\}$ is a continuous function of $x$.
Also, for a fixed $x$, $P_X\times P_X\{(y_1,y_2): d(x,\partial A(y_1,y_2))\leq 3\delta\}\to 0$, again by the dominated convergence theorem.
	
	By Lemma 1, 
	\begin{multline*}
		\sup_{x \in K}  P_X\times P_X\Big[\{y=(y_1,y_2): \omega_{f_x}\{B_\rho(y,\delta)\}\geq \epsilon\}\Big]\leq \\
		\sup_{x\in K}  P_X\times P_X\{(y_1,y_2): d(x,\partial A(y_1,y_2))\leq 3\delta\}
	\end{multline*}
	
	\noindent which converges to $0$ as $\delta\to 0$.
\QEDB

\textit{Proof of Theorem \ref{Rd}}

	Let $\epsilon>0$ and suppose that $R$ is large enough so that $P_X(\overline{B(0,R)}^c)<\epsilon$, and take $K=\overline{B(0,R)}$.
	We will apply the previously mentioned theorems of Billingsley and Tops{\o}e to the set of functions 
	$$\mathcal{F}=\{f_x(z,y)=\mathbb{I}_{A(z,y)}(x):x\in \mathbb{R}^d\}.$$ 
	We start by splitting into two terms 
	\begin{multline*}
		\sup_{x}  P_X\times P_X\Big[\{y=(y_1,y_2): \omega_{f_x}\{B_\rho(y,\delta)\}\geq \epsilon\}\Big]\leq \\ \sup_{x\in K}P_X\times P_X\Big[\{y=(y_1,y_2): \omega_{f_x}\{B_\rho(y,\delta)\}\geq \epsilon\}\Big]+\\
		\sup_{x\in K^c}P_X\times P_X\Big[\{y=(y_1,y_2): \omega_{f_x}\{B_\rho(y,\delta)\}\geq \epsilon\}\Big]=I_1+I_2.
	\end{multline*}
	
	By Lemma 1,  $I_1\leq \sup_{x\in K}  P_X\times P_X\{(y_1,y_2): d(x,\partial A(y_1,y_2))\leq 3\delta\},$ 	which converges to $0$ as $\delta\to 0$.
Regarding $I_2$, we bound
	\begin{multline*}
		P_X\times P_X\Big[\{y=(y_1,y_2): \omega_{f_x}\{B_\rho(y,\delta)\}\geq \epsilon\}\Big]\leq \\
		P_X\times P_X\Big[\{y=(y_1,y_2):y_1,y_2\in  K  \text{ and } \omega_{f_x}\{B_\rho(y,\delta)\}\geq \epsilon\}\Big]\\
		+P_X\times P_X\Big[\{y=(y_1,y_2):y_1\in K^c\text{ or } y_2\in K^c \text{ and }\omega_{f_x}\{B_\rho(y,\delta)\}\geq \epsilon\}\Big].
	\end{multline*}
	
	The second term is bounded from above by $2P_X(K^c)<2\epsilon$.
	Lastly, to tackle the first term, by Lemma 1, it only remains to  prove that
	$$\lim_{\delta\to 0}\sup_{x\in K^c} P_X\times P_X\Big[\{y=(y_1,y_2):y_1,y_2\in  K \text{ and } d(x,\partial A(y_1,y_2))<3\delta\Big]=0.$$
	
	Lastly, let $\delta_0$ be small enough such that $P_X((\partial K)^{3\delta})<\epsilon$ for all $\delta<\delta_0$.
If $y_1,y_2\in B(0,R-3\delta)$, then $d(x,A(y_1,y_2))\geq 3\delta$.
Therefore,
	
	\begin{multline}
		P_X\times P_X\Big[\{y=(y_1,y_2):y_1,y_2\in  K \text{ and } d(x,\partial A(y_1,y_2))<3\delta\Big]\leq \\ 
		P_X(y_1\in (\partial K)^{3\delta} \text{ or } y_2\in (\partial K)^{3\delta})<2\epsilon.
	\end{multline}
	
\QEDB

\

\textit{Proof of Theorem \ref{prop3}}
\\
	By assumption, we have that $\mathcal F$  has a finite $VC$-dimension.
This implies that the assumptions in Proposition $10$ in  \cite{gine1996} hold for an order $2$ $U$-statistic and the asymptotic distribution is derived from  Theorem 4.10 in  \cite{arcones1993}.
\QEDB

\

\textit{Proof of Theorem \ref{levels}}

Let $K$ be any compact set. Observe that $\{\textrm{LD}\geq \lambda^{-}\}\cap K$ is nonempty and compact for any $0<\lambda^{-}<\lambda$.
Following the same ideas as those used in the proof of Theorem 1 in \cite{cuevas2006}, one can derive that for all $\epsilon>0$
	$$\partial \{\textrm{LD}\geq \lambda\}\cap K \subset B(\partial\{\widehat{\textrm{LD}}_n\geq\lambda\}\cap K,\epsilon).$$
	The proof of
	$$\partial\{\widehat{\textrm{LD}}_n\geq \lambda\}\cap K\subset B(\partial \{\textrm{LD}\geq \lambda\}\cap K,\epsilon)$$	
	is slightly different from the analogous inclusion in \cite{cuevas2006}.
If we proceed by contradiction,   there exists an $x_n\in \{\widehat{\textrm{LD}}_n\geq \lambda\}\cap K$, but 
	\begin{equation} \label{cont1}
		d(x_n,\partial\{\textrm{LD}\geq \lambda\}\cap K)>\epsilon.
	\end{equation} 
	Observe that if $x_n\in  \partial \{\widehat{\textrm{LD}}_n\geq \lambda\}\cap K$, then $x_n\in \partial A(X_i,X_j)$ for some $X_i\neq X_j$.
	If $x_n$ is in the boundary of two or more sets $A(X_i,X_j)$,  we can take $y_n$ such that $d(y_n,x_n)<\epsilon/2^n$, $y_n\in \partial \{\widehat{\textrm{LD}}_n(x)\geq \lambda\}$,  $y_n$ in the boundary of only one $\partial A(X_i,X_j)$ and  $|\widehat{\textrm{LD}}_n(y_n)-\lambda| <\binom{n}{2}^{-1}$.
If $x_n$ is in the boundary of only one $A(X_i,X_j)$, we choose $x_n=y_n$.
In this case we clearly have also have $|\widehat{\textrm{LD}}_n(y_n)-\lambda| <\binom{n}{2}^{-1}$

	Since $K$ is compact, there exists an $x\in K$ such that $x_n,y_n\to x$ (by considering a subsequence  if necessary).
Then $|\textrm{LD}(x)-\lambda|\leq |\textrm{LD}(x)-\textrm{LD}(y_n)| +|\textrm{LD}(y_n)-\widehat{\textrm{LD}}(y_n)|+|\widehat{\textrm{LD}}(y_n)-\lambda|$.
	The first terms converge to $0$ by the continuity of $\textrm{LD}$ at $x$.
 The second term converges to $0$ a.s., because the set $K'=\{\cup_n \{y_n\}\}\cup \{x\}$ is compact and by Theorem \ref{convunif}, $\sup_{z\in K'} |\widehat{\textrm{LD}}_n(z)-\textrm{LD}(z)|\to 0$.
The last term is bounded from above by $\binom{n}{2}^{-1}\to 0$ as $n\to \infty$.
Then $\textrm{LD}(x)=\lambda$, which implies that $d(x_n,\partial \{\textrm{LD}\geq \lambda\}\cap K)\leq d(x_n,x)\to 0$, which contradicts \eqref{cont1}.
\QEDB

\bibliography{lens}

\end{document}